\documentclass[12pt,a4paper]{amsart}
\usepackage{amssymb}

\newcommand{\B}{\ensuremath{\mathcal{B}}}

\newcommand{\M}{\ensuremath{\mathcal{M}}}

\newcommand{\CC}{\ensuremath{\mathbb{C}}} 
\newcommand{\FF}{\ensuremath{\mathbb{F}}}

\newcommand{\PP}{\ensuremath{\mathbb{P}}} 
\newcommand{\pp}{\ensuremath{\mathbb{P}}}

\newcommand{\ZZ}{\ensuremath{\mathbb{Z}}} 

\newcommand{\ra}{\ensuremath{\rightarrow}} 
\newcommand{\QED}{\hspace*{\fill}\qed} 
\def\hol{{\mathcal{O}}}

\newcommand{\Proof}{{\it Proof. }}

\newtheorem{teo}{Theorem}[section] 
\newtheorem{df}[teo]{Definition}
\newtheorem{lem}[teo]{Lemma} 
\newtheorem{question}[teo]{Question}

\newtheorem{oss}[teo]{Remark} 
\newtheorem{prop}[teo]{Proposition}

\newtheorem*{nota}{Notation}

\title{The moduli space of surfaces with $K^2=6$ and $p_g=4$.}

\author{I. C. Bauer}
\author{F. Catanese}
\author{R. Pignatelli}
\begin{document}

\thanks{ The
research of the  authors was performed in the realm  of the
  DFG SCHWERPUNKT "Globale Methode in der komplexen Geometrie", and of the
EAGER EEC Project. The third author was  supported by the Schwerpunkt and by P.R.I.N. 2002
"Geometria delle varietˆ algebriche" of  M.I.U.R. and is a member of G.N.S.A.G.A. of
I.N.d.A.M. }

\maketitle
\tableofcontents

\section*{Introduction}

The motivation for our work stems from the questions posed by F. Enriques
in Chapter VIII of his book "Le superficie algebriche" (\cite{enr}) about
surfaces of general type with $p_g = 4$ and their moduli. 

For these, $K^2 \geq 4$, and the cases  $K^2 = 4,5$ were completely
classified by Enriques(cf. \cite{enr}). Enriques also discussed at length the case
$K^2=6$, which was later completely classified by Horikawa in \cite{horIII}.

The first question posed by Enriques was the following: for which value 
of $K^2$ does there exist a surface with $p_g = 4$ and birational canonical map?
This existence question, posed by Enriques for 
$K^2 \geq 8$, was later solved by virtue of the contributions of several authors, and we now
know that such surfaces exist for  $ 7 \leq K^2 \leq
32$, (cf. e.g. \cite{cil}, \cite{warsaw}).

The answer to the question concerning classification and moduli is much
harder, and a complete classification has been achieved up to now only for 
$K^2 \leq7$, see for instance  the monograph (\cite{7}) for the case $K^2 = 7$.

The challenging open problem for $K^2 = 6,7$ is to completely understand the structure of
the moduli space, i.e., to determine the incidence correspondence of the several
locally closed strata which are described in the classification.

Horikawa  showed in \cite{Q} that the moduli space for $K^2=5$ is connected,
with two irreducible components of dimension $40$ meeting along a divisor.

 He showed later  (\cite{horIII}) that for  $K^2=6$ there are exactly four
irreducible components and at most three connected components . 

He did so by first listing all possibilities for the canonical
map, dividing thus the moduli space of surfaces with $K^2=6$ into $11$  nonempty locally
closed strata, and then analysing some of their local deformations.

More precisely, Horikawa named the $11$ strata 
$I_a$, $I_b$, $II$, $III_a$, $III_b$, $IV_{a_1}$, $IV_{a_2}$,
$IV_{b_1}$,  $IV_{b_2}$, $V_1$, $V_2$ 
(see \cite{horIII} or (\ref{strata}) below for precise definitions of each stratum). 
According to Horikawa's notation we define 
\begin{nota}
Let $A$ and $B$ be two of the above introduced strata. 
The notation ``$A \rightarrow B$'' means that $B$ intersects the closure of $A$,
 i.e., there is a deformation of a surface of type $B$ to surfaces of type $A$
(it suffices to have a flat family over a small disk $\Delta_{\varepsilon} \subset
\CC$  whose central fibre is of type $B$ and whose general fibre is of type $A$).
\end{nota}
With this notation Horikawa summarized his results in the following
diagram
$$
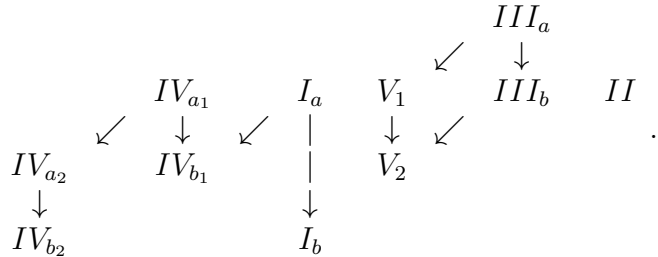
\begin{array}{ccccccccccc}
&&&&&&&&III_a&&\\
&&&&&&&\swarrow&\downarrow&&\\
&&IV_{a_1}&&I_a&&V_1&&III_b&&II\\
&\swarrow&\downarrow&\swarrow&|&&\downarrow&\swarrow&&&\\
IV_{a_2}&&IV_{b_1}&&|&&V_2&&&&\\
\downarrow&&&&\downarrow&&&&&&\\
IV_{b_2}&&&&I_b&&&&&&
\end{array}.
$$

 Our main result is:

\begin{teo}\label{Degree6}
Consider the moduli space of surfaces with $p_g=4, K^2 = 6$.
Then it has at most two connected components.

In particular, $ II \ra III_b$, i.e., there is a deformation of surfaces of type
($III_b$) to surfaces of type (II).

\end{teo}

We would moreover like to pose the following

\begin{question}
It  the above moduli space  
(for  $ p_g=4, K^2=6$) disconnected?
\end{question}

A possible reason for this could be that the surfaces of all irreducible components
degenerate to surfaces with a genus two pencil, but in one case the
braid monodromy is transitive, in the other case it is not.

The new idea that we exploit is the following: 
the canonical models  $X$ of surfaces of Type
(II) are exactly the hypersurfaces of degree
$9$ in the weighted projective space $\pp(1,1,2,3)$ (with rational double
points  as singularities). This remark was used in  \cite{nonred}) to give
a new explanation of the result of Horikawa that the moduli space is non reduced on
the open set (II), and it implies that the canonical divisor is 2-divisible as a
Weil divisor on
$X$.

Similarly occurs for type ($III_b$), so for both type of surfaces we have a
semicanonical ring $\B$, and what we do is to find a flat family of
deformations of the semicanonical ring. 

How to do this?
The ring $\B$ is a Gorenstein ring, of codimension $1$ in case (II),
of codimension $4$ in case ($III_b$), where $X$ is embedded in  
$\pp(1,1,2,3,4,5,6)$.

In order to describe the semicanonical ring and its deformations in case
($III_b$), we use, as in
\cite{bcp}, the format of $4 \times 4$ Pfaffians of 
extrasymmetric antisymmetric $ 6
\times 6$ matrices. 

This format applies to a codimension $2$ subvariety of the stratum of surfaces of type
($III_b$): these  have a  pencil of hyperelliptic curves of genus $3$, and one has
to lift  (cf. \cite{hyp}) this graded
ring of dimension $1$ to the semicanonical ring of the surface.

The deformation trick is similar to the one used in \cite{bcp} for the 
canonical ring: filling entries of homogeneous degree $0$ in the matrix
with parameters (and respecting then the symmetries). When these parameters are non
zero, three of the given Pfaffians allow to eliminate the $3$ variables of
respective weights
$4,5,6$. We obtain then a semicanonical ring of type (II).

We want now to briefly discuss the cited method of
 extrasymmetric antisymmetric $ 6 \times 6$
matrices. 

 The main point here is the
lack of a structure theorem for Gorenstein subvarieties of codimension $4$
(for codimension $3$ we have the celebrated theorem of Buchsbaum and
Eisenbud \cite{b-e}). 

Several  explicit formats were proposed by Dicks, Reid and Papadakis.

The geometric roots (cf. \cite{graded}) for the format we use here
lie in the fact that the Segre product
$\pp^2 \times \pp^2$ is embedded in $\pp^8$ as the variety of $3 \times 3$
matrices $A$ of rank $1$, hence defined there by $9$ quadratic equations,
admitting $16$ syzygies. 

If however one writes $ A = B + C$, with $C$ symmetric, $B$
antisymmetric, then one can form the   antisymmetric $6 \times 6$-matrix:

$$
D=\begin{pmatrix}
B& C  \\
- C& B  \\
\end{pmatrix}.
$$

The matrix $D$ has an extrasymmetry from which follows indeed that the $15$
($4 \times 4$) Pfaffians of
$D$ are not  independent, but  exactly reduce to the  above $9$
quadratic equations.

Using a flat family of deformations of the above subvariety, and
interpreting the entries of the matrix as indeterminates to be
specialized, one  obtains an easy construction of Gorenstein subvarieties of
codimension $4$.

We refer to \cite{graded} for a thorough discussion of the problem of
understanding Gorenstein rings in codimension $4$. Our result shows that the
"moduli space" of such rings could be rather complicated, since this format does not
apply for a general surface of type ($III_b$), and moreover since
we obtain a deformation from codimension $4$ to codimension $1$, but we observe that
one cannot pass through all the lower codimensions.

In the next section we shall determine which of the eleven Horikawa classes yield canonical
models $X$ where the canonical divisor is 2-divisible as a Weil divisor
(we use here the classical notation $\equiv$ for linear equivalence).

\section{Surfaces with 2-divisible canonical divisor}

Let $S$ be the minimal model of a surface with $ K_S^2=6, p_g=4$, and let $X$ be its
canonical model (obtained by contracting the curves $C$ with  $K_S \cdot C = 0$).

In this section we shall classify the closed set of the moduli space given by the surfaces
for which 

(***) {\bf there exists a Weil divisor $ L $ on $ X $ such that $ 2 L
\equiv K_X $ and $h^0 (X, L) \geq 2$.}

Observe that the above hypothesis immediately implies that the image of the canonical map
$\phi_{K_X}$ is a quadric cone and that $h^0 (X, L) = 2$.

Recall now
that, if $|K_S| $ has a nonempty fixed part $\Phi$ and we write
$|K_S| = |M| +
\Phi$, then, since
$K_S$ is nef,
$ 6 = K_S^2 = K_S \cdot \Phi + M \cdot \Phi + M^2 \geq M \cdot \Phi + M^2 \geq 2 + 4 $,
the last inequality following from the 2-connectedness of canonical divisors, and by the
fact that the canonical system is not a pencil(as shown in \cite{horIII}).

Whence, if the fixed part $\Phi \neq \emptyset$ we have
$$ M^2 = 4,  M \cdot \Phi = 2,  K_S \cdot \Phi = 0.$$ 

Therefore we conclude that $K_X$ has no fixed part on the canonical model $X$, and
only the following cases a priori occur:

\begin{itemize}
\item
(0) $|K_X|$ is base point free
\item
(2) $|K_X|$ has 2 smooth base points
\item
$|K_X|$ has a smooth base point plus possibly an infinitely near one (these 
two cases are however
shown by Horikawa not to occur)
\item
(1) $|K_X|$ has a singular base point $p$ on $X$.
\end{itemize}

The last case is the only one where there is a fixed part $\Phi$ on the minimal model.
Since to the  point $ p \in Sing(X)$ corresponds the fundamental cycle $Z$ (pull back of
the maximal ideal $\M_p \subset \hol_{X,p}$), we have $\Phi \geq Z$.
$Z$ has the property that for each divisor $\Psi$ which lies in the inverse image of
$p$ one has $Z \cdot \Psi \leq 0$

Write $ \Phi = Z + \Psi$, and assume $\Psi >0$. Observe then that $M \cdot \Phi = 2,  K_S
\cdot \Phi = 0$ implies $\Phi^2 = -2$. Then $ -2 = \Phi^2 = Z^2 + 2 Z \cdot \Psi +
\Psi^2 
\leq -2 + 0 -2 = -4$ is a contradiction, whence $\Phi$ equals the fundamental cycle.

Denote by $Q$ the image of the canonical map and assume that it is a quadric cone.

 Blowing
up the base points of $\phi_{K_X}$ we obtain a morphism $f^0 : X^0 \ra Q$, and let us
observe that the irreducible exceptional curves of $ X^0 \ra X$ are  at most 2.

In any case the blow-up formula yields

\begin{itemize}
\item
$K_{X^0} = \pi^*(K_X)$ in cases (0) and (1)
\item
in case (2) there are two (-1)-curves $E_1, E_2$ such that $K_{X^0} -E_1 - E_2 =
\pi^*(K_X)$.
\end{itemize} 

\begin{oss}
In case (0) the canonical map $\Phi_X$ is a finite morphism, whence (***) holds
if and only if $ Im (\Phi_X) $ is a quadric cone, i.e., if and only if $\Phi_X$
has degree 3 (as it is easy to see, cf. Lemma 4.1 of \cite{horIII}).

In the other two cases we have $ deg (\Phi_X) = 2$.
\end{oss}

\begin{oss}

Assume that we are in cases (2), (1) and that the canonical image is a quadric cone $Q$.

Let $L'$, $L^0$, $L''$ be the respective proper transforms on $X$, resp. $X^0$, resp.
$S$ of a general line $ l \subset Q$. 
 
Observe that the canonical divisor $K_X$ is the pull back of a hyperplane divisor
on $Q$, and its pull back to $X^0$ has as movable part the pull back $H$ of 
a hyperplane divisor on $Q$. 

Whence  $ H \equiv 2 L^0 + W$, where
$W$ is the effective divisor on
$X^0$ corresponding to the inverse image of the vertex $v \in Q$.

More precisely, $W$ is the fixed part of $ |H - 2 L^0|$. 

 Assume that there is a Weil divisor $L$ satisfying (***): then one immediately sees
that $h^0(X,L)=2$ and there is an effective divisor $E'$ on $X$ with $ L \equiv L' + E'$,
and where $E'$ is the fixed part of the pencil $|L|$. Then the fixed part of $K_X - 2L'$
equals $2 E'$, whence $ \pi_*(W) = 2 E'$.

It follows that (the linear equivalence class of) $K_X \equiv 2 L' +   \pi_*(W)$ is
2-divisible as a Weil divisor if and only if $\pi_*(W)$ is 2-divisible as an effective
divisor.

There are two possibilities for this: ${f^0}^{-1}(v)$ has dimension 0, or , in case where
${f^0}^{-1}(v)$ has dimension 1, we have $\pi_*{f^0}^{-1}(v) = 2 E'$. In this last case, we
have that $f^0$ factors through the Segre-Hirzebruch surface  $\FF_2$.
\end{oss}

We recall now the case subdivision given by Horikawa (\cite{horIII})

\begin{df} \label{strata}
Assume we are in case (0) ($K_X$ has no base points): then we have   type ($I_a$)
if $\phi_{K_X}$ has degree 1, type (Ib) if $\phi_{K_X}$ has degree 2,
type (II) if $\phi_{K_X}$ has degree 3. 

Case (III) is the case where  $\phi_{K_X}$ has degree 2, but there is
no genus $2$ pencil on $X$. There are two subcases: ($III_a$), where the canonical
image is a smooth quadric and we have two smooth base points, and ($III_b$), where 
 the canonical image
is a  quadric cone and we have one singular base point.

The two cases of type (IVa-1), (IVa-2) have a smooth quadric as canonical image,
the two cases of type (IVb-1), (IVb-2) have a  quadric cone as canonical image, 
$X^0$ is a double cover of $\FF_2$ , but the section at 
infinity is not part of the branch locus. 

Surfaces of type (IV) and (V) have all
a genus 2 pencil. Surfaces of type (V-1) (V-2) have a quadric cone as canonical image,
$X^0$ is a double cover of $\FF_2$ , and the section at 
infinity is  part of the branch locus. For type (V-1) we are in case (2), for type (V-2) 
we are in case (1).

\end{df}

\begin{prop}
The canonical model $X$ of a surface with $K^2 = 6, p_g=4$ satisfies condition (***) 
(there exists a Weil divisor $ L$ on $ X $ such that $ 2 L
\equiv K_X $ and $h^0 (X, L) \geq 2$) if and only if it is of one of the following
types: (II), ($III_b$), (V-1) or (V-2).
\end{prop}

\proof
Since the canonical image must be a quadric cone $Q$, cases (Ia), (Ib), ($III_a$),
(IVa-1), (IVa-2) are immediately excluded.

Assume that (***) holds for some surface of type  (IVb-1) or (IVb-2). We
know that the section at infinity
$\Delta_{\infty}$ is not part of the branch locus on $\FF_2$.

Then, by our previous remark, the inverse image 
$W_{\infty}$ of $\Delta_{\infty}$ on $X^0$ must be contracted by $\pi : X^0 \ra X$.

In case (2), we observe that $|H|$ is base point free and $ H \cdot E_i = 1$, whence
$E_i$ maps to a line and not to the vertex.

In case (1), again $|H|$ is base point free, and by \cite{horIII}. theorem 6.2 the base
point of $|K_X|$ is an ordinary double point. Let $F$ be the corresponding $-2$ curve: then
$H \cdot F = 2$ and we have again a contradiction.

We have already seen that (***) holds for type (II). 

In case ($III_b$) the section $\Delta_{\infty}$ (cf. Theorem 5.2, ibidem) is
isolated in the branch locus, whence its set theoretic inverse image is a smooth
$-1$-curve 
 $W_{\infty}$, thus  ${f^0}^{-1}(v)$ has dimension 0.

In cases (V-1), (V-2) again one shows that ${f^0}^{-1}(v)$ has dimension 0, using Theorems
6.1, resp. 6.2, ibidem.

\QED

\section{Surfaces of type $II$ and $III_b$.}

In this section we want to concentrate on the above classes of surfaces, the ones
for which (***) holds, but there is no genus $2$ pencil on $S$.

The following lemma shows that the classes (II) and ($IIIb$) are exactly those 
for which the pencil $|L|$ has no fixed part.

\begin{lem}
Assume  that (***) holds and  write $ K_S \equiv 2 \Lambda + Z$ where $K_S \cdot Z =
0$. Then the pencil $|L|$ is without fixed part  (i.e., $L' \equiv L$, equivalently
$\Lambda \equiv L''$)
if and only if there is no genus $2$ pencil on $X$.
In this case,  the general element in $|L''|$ is smooth
irreducible  of  genus $g(L'')= 3$ and we have: $(L'')^2 = 1, L'' \cdot Z
= 1$. It follows also that $Z^2= -2$ and thus $Z$ is the fundamental cycle
 of a singular point $P_1$ of $X$.
\end{lem}

\Proof
Write $|\Lambda| = |L''| + E''$ with $E'' > 0 $.
 Since $K_S \equiv 2 L'' + 2 E'' + Z$,
we have $ L'' \cdot K_S + E'' \cdot K_S = 3$, and moreover 
$ L'' \cdot K_S > 0,  E'' \cdot K_S > 0.$

It is impossible that $ L'' \cdot K_S = 1$, since then $ {L'' }^2$ is odd,
and $ L'' \cdot K_S \geq 2$ ($L''$ is nef), a contradiction.

Thus $ L'' \cdot K_S = 2$ and $ {L''} ^2 = 0$, thus we have a genus 2 pencil.

Conversely, the curves of a genus 2 pencil map to the lines of the quadric cone $Q$,
but if  $|L|$ is without fixed part then $E''= 0$, and we claim that $L'' = \Lambda$
is a pencil of genus 3 curves. 

In fact, $L'' K_S = 3 = 2 (L'')^2 +  L' Z$, thus  $ L'' Z$ is odd. Since $ L'' Z$
is non negative  and odd, while   $L''
K_S = 3 $ implies that 
$(L'')^2 $ is also odd,  the only possibility is that
$(L'')^2  = 1, L'' Z = 1$. Whence, $ p(L'') = 3$.

In particular,  $|L''|$ has a unique smooth base point $P_0$ and the general 
curve in $|L''|$ is smooth by Bertini's theorem.

Since $L'' Z = 1$, follows that $Z^2 = -2$, and since $Z$ is the only 
 divisor in $ 2 L'' + Z \equiv K_S$ exceptional for $ S \ra X$, follows that
$Z$ is a fundamental cycle.

\qed

\begin{df} \label{semi}
Assume that (***) holds.Then the {\bf semicanonical ring of $X$} is the
graded ring
$$ \B := \oplus_{m=0}^{\infty}  \B_m := \oplus_{m=0}^{\infty} H^0(X, mL).$$
\end{df}

\begin{oss}
Obviously, 
$$ \B_{2m} \cong  H^0(X, mK_X) \cong  H^0(S, mK_S), \B_{2m+1} 
\cong  H^0(S, mK_S + L'').$$
\end{oss}

\begin{lem}\label{odd}
Assume that $|L|$ is without fixed part on $X$. 
Then the  sequence
$$ 0 \rightarrow H^0(S, \hol(m K_S)) \rightarrow H^0(S, \hol(m K_S+L'') )
\rightarrow  H^0(\omega_{L''}((m-1)K_S)) \rightarrow 0$$
is exact.
Moreover, $dim \B_{2m} = 5 + 3 m (m-1)$ for $m \geq 2$, $dim \B_{2m+1} = 7 + 3
(m+1) (m-1)$ for $ m\geq 1$. 
\end{lem}
\Proof
We have an exact sequence of sheaves given by the adjunction formula and moreover
$h^1(\hol(mK_S)) = 0$ for each $m$ since $q(S)=0$. The rest follows from the previous
remark and from the fact that
$ K_S \cdot L'' = 3$ and a general $L''$ is smooth of genus $3$.

\QED

\begin{lem}\label{3L}
Consider a basis $\{ x_0, x_1 \} $ of $H^0 (X, L) $  and  pick $y_2 \in
H^0 (X, 2L)$ in order to complete 
 $\{ x_0^2, x_0 x_1, x_1^2 \} $ to a basis of $H^0 (X, 2L)$.
Then $\{ x_0^3, x_0^2 x_1, x_0 x_1^2  , x_1^3, x_0 y_2 , x_1 y_2\}$ 
are linearly independent and there exists an element $z_3 \in
H^0 (X, 3L)$ completing the previous set to a basis
of $H^0 (X, 3L)$.
\end{lem}

\Proof

Otherwise, w.l.o.g. we may assume that we have a relation of the form
$$ x_0 y_2 = G_3(x_0, x_1)$$
where $G$ is homogeneous of degree $3$ and divisible by $x_1$.

But then, setting $y_{00} := x_0^2, y_{11} := x_1^2, y_{01} := x_0 x_1$,
we have $ G_3(x_0, x_1) = x_1 L(y_{00}, y_{11}, y_{01},
y_2)$ where $L$ is a linear form. Whence, the canonical image  of
$X$ satisfies $3$ quadratic equations 
$$y_{00}y_{11}=  y_{01}^2, y_{00}y_2 = y_{01} L(y), y_{01}y_2 = y_{11} L(y). 
$$

But it is shown by Horikawa (\cite{horIII}) that $|K_S|$ is not a pencil.

There remains   only to observe that $h^0 (X, 3L) = 7$, by the previous 
lemma.

\qed

\begin{lem} \label{II}
Assume that $|L|$ is without fixed part on $X$.
The following three conditions are equivalent:
\begin{itemize}
\item
$X$ is of type (II), i.e., the canonical map has degree 3   
\item
 the general curve in $|L''|$ is non hyperelliptic
\item
 the sections $ x_0 ,x_1 , y_2,
z_3 $ provide a morphism $\psi$ to the 
$3$-dimensional  weighted projective space $\PP (1,1,2,3)$ which is
birational onto a hypersurface $\Sigma$ of degree $9$.
\end{itemize}

 In case (IIIb) we have:
\begin{itemize}
\item
 every  curve in $|L''|$ is hyperelliptic,
 \item 
 the base point $P_0$ of $|L''|$ is a Weierstra\ss\ point of 
 every curve $L''$
 \item
 $Z$ is the fixed part of $K_S$
 \item 
 the canonical map is a double covering of a quadric cone.
\end{itemize}
\end{lem}

\Proof
Since $|L''|$ has a simple base point $P_0$ the general curve in the pencil
is smooth of genus $3$ and its canonical system has no base points. Whence
$|3 L'' + Z| = | K_S + L''|$ has no base points on the general curve $L''$, and in
particular $P_0$ is not a base point; thus the base points are contained in $Z$.
 Whence on $X$ the only possible base
 point (where $x_0 = x_1 = z_3 = 0$) is  $P_1$. 
 
 We observe that there is a base point of $\psi$ if and only if $y_2$
vanishes on $P_1$, and this case will be denoted by (*). 
  
  In case (*) $P_1$ is a point where $x_0 = x_1 = y_2 = 0$, whence
$Z$ is  in the fixed part of 
  the canonical system, but $|K_S - Z| = |2L''| $ has no base point since $(K_S
- Z)^2 = 4$ and its base locus is  $\subset \{P_0 \}$. Hence, $|K_S - Z|$
yields a  morphism $f : X \rightarrow \PP^3 $ which is a double cover of the
quadric cone, and $P_1$ is the  unique point where $x_0 = x_1 = y_2 = 0$,
and $Z$ is exactly the fixed part of 
  the canonical system.
 
 In   general, by lemma \ref{3L} follows that $\psi$ is birational if and only
if the general curve in $|L'|$ is non -hyperelliptic, and in any case the
degree of $\psi$ is at most $2$. 

If we are not in case  (*), we simply observe that $\psi$ is  then 
a morphism
and that $ deg (\psi) deg (\Sigma) = 9 = 6 ( 3/2)$, but we
have already seen that if the map is not birational, its degree is $2$.

Whence, it follows that  $\psi$ is a  birational morphism (and obviously
then $deg (\Sigma) = 9$).

 In case (*) every curve in $|L'|$ is a double cover of a line, whence 
 all the curves in
$|L'|$ are hyperelliptic. Since the point $P_0$ is invariant by the 
involution of $S$ yielding the hyperelliptic involution on every curve in 
 $|L'|$, it follows that $P_0$ is a Weierstra\ss\ point.
 
More precisely, the exact sequence

$$  0 \rightarrow \hol_S (L'' + Z)  \rightarrow \hol_S (K_S) 
\rightarrow \hol_{L''} (K_S) \rightarrow 0$$
and the remark that  $ \hol_{L''} (K_S) = \omega_{L''} (-P_0)$
shows that $|K_S|$ has $P_0$ as a base point on $L''$ if and only if $L''$
is hyperelliptic and $P_0$ is a Weierstra\ss\ point.

\qed

Horikawa gave a very concrete description of  surfaces of type (IIIb).

\begin{teo}[\bf 5.2 in \cite{horIII}] \label{IIIb}
Let $S$ be a surface of type $III_b$. Then $S$ is birationally
equivalent to a double covering of $\FF_2$ whose branch locus $B$
consists of the negative section $\Delta_{\infty}$ and of $B_0 \in |7\Delta_{\infty}+14
\Gamma|$ which has a quadruple point $x$ and a $(3,3)$- point at $y
$ such that $x$ and $y$ belong to the same fibre $\Gamma_0 \in |\Gamma |$. Moreover, $y$ may be infinitely near to $x$.
\end{teo}

\section{On rings associated to curves of genus $3$}

To compute the semi-canonical ring we use ideas related to the hyperplane section
principle introduced by Miles Reid (cf. page 218 of
\cite{hyp}). 

{\bf The hyperplane section principle:} Let $\B$ be a graded ring, and
$x_0 \in \B$ a homogeneous non-zero divisor of degree $\deg{x_0} >0$; set
$\overline{\B}=\B/(x_0)$. The hyperplane section principle says that quite
generally, the generators, relations and syzygies of $\B$ reduce mod
$x_0$ to those of $\overline{\B}$.

\begin{prop}
Let $\B$ be the semicanonical ring of $X$ and  fix  an element $x_0$ of degree $1$
in
$\B$ whose divisor yields a smooth curve $C \in |L''|$.   
Then the quotient ring  $\overline{\B} = \B / (x_0)$ satisfies
$$\overline{\B}_{2m+1} = H^0(\omega_{L''}((m-1)K_S)), \overline{\B}_{2m} =
H^0(\hol_{L''}(m K_S)). $$
\end{prop}
\Proof
The first assertion follows immediately from  lemma \ref{odd}.
The second assertion is immediately verified for $m=1$, while, for
$ m \geq 2$, it follows from the exact sequence
$$ 0 \rightarrow H^0(S, \hol(m K_S - L'')) \rightarrow H^0(S, \hol(m K_S) )
\rightarrow  H^0(\hol_{L''}((m)K_S)) \rightarrow $$
$$ \ra H^1(S, \hol(m K_S - L''))
\ra 0$$ and the vanishing of 
$$h^1(S, \hol(m K_S - L'')) = h^1(S, \hol(-[(m-2) K_S +
L'' + Z])$$.

Here we use Serre duality plus the fact that, on a regular surface ($q=0$)
$H^1(S, \hol(-D)) =0$ if the divisor $D$ is effective and numerically connected
(cf. \cite{CM}). That $(m-2) K_S +
L'' + Z$ is numerically connected 
can be proved directly, but follows more easily for case (IIIb) when we observe
that by Theorem 5.2 of \cite{horIII}, $Z$ is an irreducible $-2$-curve.

\qed

Let us compute the ring $ \overline{\B}$  for surfaces of type $III_b$.

We see immediately  from the previous proposition, and from lemma \ref{II}
that the quotient ring $\overline{\B}$ is isomorphic to a ring
of the type described in the following

\begin{df}
Let $C$ be a smooth hyperelliptic curve of genus $3$, $p\in C$ be a Weierstra\ss\ 
point, and $X $ a section of $H^0(\hol_C(p)$ with $div(X)=p$. 

Consider  the Itaka ring $R(C,p): = \oplus_{d \geq 0} H^0(\hol_C(p)^{\otimes n})$
and  
define $R(C,\frac32 p)$ as the graded subring with
$$
\left\{
\begin{array}{lll}
R(C,\frac32 p)_{2d}&:=&R(C,p)_{3d}\\
R(C,\frac32 p)_{2d+1}&:=&R(C,p)_{3d+1} 
\end{array},
\right .
$$
and with product defined,
for homogeneous elements, as $ab= a \otimes b$, resp. $a \otimes b
\otimes X$ according to the parity of the product of the degrees of $a$ and $b$
(even, resp.  odd). 
\end{df}

So, the ring $ \overline{\B}$  for surfaces of type $III_b$ being a subring of
$R(C,p)$, we need first to describe the latter.

The ring of a Weierstra\ss\ point of a smooth hyperelliptic curve is
well known in every genus: for the convenience of the reader we state
and prove here the result in the case of genus $3$.

\begin{lem}\label{RCp}
Let $C$ be a hyperelliptic curve of genus $3$, $p \in C$ a
Weierstra\ss\ point: then $R(C,p)\cong \CC[X,Y,T]/(T^2-P_{14}(X,Y))$ where
$\deg(X,Y,T)=(1,2,7)$ and $P_{14}$ is homogeneous of degree $14$.
\end{lem}

\Proof
Let $X $ be a section of $H^0(\hol_C(p))$ with $div(X)=p$: the section $X$ is
antiinvariant and the divisor
$p$ is invariant under the hyperelliptic involution $\sigma$, such that $ \phi: C
\ra C / \sigma \cong\PP^1$ is branched on a divisor $B$  of degree $8$.

The morphism $\phi$ is given by a basis of $H^0(\hol_C(2p))$, for instance let us
take $X^2$ and a new element $Y$.

We have $\phi_* \hol_C = \hol_{\PP^1} \oplus Z\hol_{\PP^1}(-4) $, where
$Z$ is an equation for the ramification divisor of $\phi$.

The even part of our ring is thus $ \oplus_{m=0}^{\infty} H^0( \phi_* \hol_C(m))
= \CC[X^2,Y] \oplus Z \CC[X^2,Y] $.

Since $Z \in H^0(\hol_C(8p))$, and $X$ divides $Z$, we may write
$$  Z = X T, T \in H^0(\hol_C(7p)).$$

Observe that $X,Z$ being antiinvariant, then $T$ is invariant. 

Consider now $H^0(\hol_C((2m+1)p)$: this space splits as the direct sum of the
$(+1)$, respectively $(-1)$-eigenspace.
By looking at the behaviour at the ramification points, we see immediately that the
sections of $H^0(\hol_C((2m+1)p)^+$ are divisible by $T$, the ones of
$H^0(\hol_C((2m+1)p)^-$ are divisible by $X$, thus the odd part of our ring
is $ T\CC[X^2,Y] \oplus X \CC[X^2,Y] $.

It follows that our ring is $ \CC[X,Y] \oplus T \CC[X,Y] $.

Its ring structure is easily obtained when we observe that $T^2$ is the pull back
of the equation of the seven remaining  branch points: 
thus we have a relation of the form $T^2= P_{14}(X,Y)$ and our claim is proven.

\QED

\begin{prop}\label{primamatrice}
Let $C$ be a hyperelliptic curve of genus $3$, $p\in C$ be a Weierstra\ss\ point. 
Then $R(C,\frac32 p)\cong \CC[x,y,z,w,v,u]/I$, where 
$\deg(x,y,z,w,v,u)=(1,2,3,4,5,6)$ and the ideal $I$ is generated by
the $4 \times 4$ pfaffians of the skew-symmetric 'extra-symmetric' matrix 
$$M= \begin{pmatrix}
  0 &0&z& v & y  & x  \\
    &0&w& u & z  & y  \\
    & &0& \tilde{P}_9 & u  & v  \\
    & & & 0 &w^2 &zw\\
    & & &   & 0  & 0\\
-sym& & &   &    & 0
\end{pmatrix},$$
where $\tilde{P}_9$ is a homogeneous polynomial of degree $9$ in the variables
$x,y,z,w$. 
\end{prop}

\Proof

It is easy to verify that the following  $6$ elements
$$\begin{array}{ccccccccc}
x&=&X&y&=&XY&z&=&Y^2\\
w&=&Y^3&v&=&T&u&=&YT
\end{array};$$
 generate $R(C,\frac32 p)$; note that
by definition the $2 \times 2$ minors of the matrix
\begin{equation}\label{4x2}
\begin{pmatrix}
x&y&z&v\\
y&z&w&u
\end{pmatrix}
\end{equation}
belong to the ideal $I$ of relations.

The other relations for these generators  come from the equation
$T^2-P_{14}$; we note that the coefficient of $Y^7$ in the polynomial $P_{14}$
of lemma \ref{RCp} cannot vanish (or the $8$ branch points of the
hyperelliptic map $\phi$  would not be distinct), and therefore we can
assume w.l.o.g. that this coefficient  be $1$. We write then
the relation as $T^2-Y^7-XP_{13}$. Let $\tilde{P}_9$ be a homogeneous
polynomial of degree $9$ in the variables $(x,y,z,w)$ such that 
$\tilde{P}_9 (X, XY, Y^2, Y^3) = P_{13}(X,Y)$ (the reader can easily check that it
exists and is uniquely determined modulo the $2 \times 2$ minors of the matrix
(\ref{4x2})).

 The relation $T^2-P_{14}$ is a relation in degree $14$ in $R(C,p)$,
and $R(C,p)_{14}$ is not contained in the subring $R(C,\frac32 p)$:
but multiplying it by suitable monomials we obtain the following relations 
\medskip

\begin{tabular}{cccccc}
$XT^2$&-&$XY^7$&-&$X^2P_{13}(X,Y)$& in degree $15$\\
$YT^2$&-&$Y^8$&-&$XYP_{13}(X,Y)$& in degree $16$\\
$Y^2T^2$&-&$Y^9$&-&$XY^2P_{13}(X,Y)$& in degree $18$,
\end{tabular} 

 which can be rewritten as polynomials in the
variables $x,y,z,w,v,u$ (uniquely determined modulo the $2 \times
2$ minors of (\ref{4x2})) as

\begin{tabular}{cccccc}
$v^2$&-&$z^2w$&-&$x\tilde{P}_{9}$& in degree $10$\\
$vu$&-&$zw^2$&-&$y\tilde{P}_{9}$& in degree $11$\\
$u^2$&-&$w^3$&-&$z\tilde{P}_{9}$& in degree $12$.
\end{tabular} 

The reader can immediately check that the ideal $I$ generated by these three
last equations and the $2 \times 2$ minors of (\ref{4x2}) coincides
with the ideal of the $4 \times 4$ pfaffians of the matrix in the
statement. 

There are no other relations since one can easily check that the
Hilbert function of $R(C, \frac32 p)$ coincides with the one
of $\CC[x,y,z,w,v,u]/I$.

\QED

\section{The family of deformations}

The hyperplane section principle gives a strategy in order to
reconstruct the ring $\B$ for every surface of type $III_b$:
taken $C$ a smooth element of the pencil $|L''|$ and $x_0$ a corresponding
section the 'hyperplane section' quotient $ \overline{\B} =  \B /(x_0)$
equals the ring $R(C,\frac32 p)$ in proposition \ref{primamatrice}.

$\B$ is obtained from $R(C,\frac32 p)$ by adding
the generator $x_0$ and deforming the $9$ equations adding suitable
multiples of $x_0$ in such a way that all the syzygies also deform.

To compute all possible 'extensions' as above is in general very difficult.
For these problems it is in
general useful to use a 'flexible format' (i.e., with free parameters) as the
one we are going to recall.

{\bf The extra-symmetric format.}
Let $A$ be a polynomial ring.

A skew-symmetric matrix $M$ is said to be
 {\it extra-symmetric} if it has the form
$$
\begin{pmatrix}
0&a&b&c& d& e\\
 &0&f&g& h& d\\
 & &0&i&pg&pc\\
 & & &0&qf&qb\\
 & & & &0&pqa\\
-sym&&&&&0
\end{pmatrix}
$$

for suitable elements $a,b,c,d,e,f,g,h,i,p,q \in A$. Then the ideal of the fifteen $4
\times 4$ pfaffians is generated by  $9$ of them; moreover, if the
entries are general enough, this pfaffian ideal  has
exactly $16$ independent syzygies, {\it which can all be 
  explicitly written as functions of the entries of the matrix} (the
computation is done in \cite{hyp} in a slightly more special case, and
it extends to this more general case).

This implies that, if we have a ring that can be written in this form,
and the ring has no further syzygies (except the $16$ we know),
deforming the entries of the matrix (preserving the symmetries)
we get at the same time a
deformation of the ideal and a deformation of the syzygies (i.e. a
{\it flat} deformation).

The first example of a ring presented through the $4 \times 4$ pfaffians of a
skew-symmetric extra-symmetric matrix was produced by  D. Dicks and
M. Reid (\cite{dicks}). This more general form appeared in
\cite{graded}; see  also \cite{bcp} for another application of it.

\

In lemma \ref{primamatrice} we wrote the ring $R(C, \frac32 p)=
\overline{\B}$ in extra-symmetric format; therefore, if we add a
variable $x_0$ in degree $1$ and lift $M$ to a matrix $N$ that is
still skew-symmetric and extra-symmetric, its pfaffians should define a
surface of type $III_b$.

Note that, since the format is not complete, the family which we
write can be (and in fact  is) smaller than the whole family of
surfaces of type $III_b$.

We add then the variable $x_0$ and rename the old variable $x$ by $x_1$.

\begin{teo}\label{famigliadidimensione36}
Consider the ring  $\CC[x_0,x_1,y,z,w,v,u]$ with variables of
respective degrees  $(1,1,2,3,4,5,6)$.

Consider a skew-symmetric extra-symmetric matrix 
$$M'=
\begin{pmatrix}
  0 &0&z& v & y  & x_1  \\
    &0&w& u &P_3 & y  \\
    & &0&P_9& u  & v  \\
    & & & 0 &wP_4&zP_4\\
    & & &   & 0  & 0 \\
-sym& & &   &    & 0
\end{pmatrix}.$$ 
where the $P_i$'s are homogeneous polynomials of degree $i$ in the first $5$
variables of the ring, and let $J$
be the ideal generated by the $4 \times 4$ pfaffians of $M$.
 
Then, for general choice of the polynomials $P_i$, $\CC[x_0,x_1,y,z,w,v,u]/J$ 
is the semi-canonical ring of a surface of type $III_b$: we obtain in this way
a family of codimension $2$ (and therefore of dimension $36$) in the
stratum of the moduli space of surfaces with $p_g=4$, $K^2=6$
corresponding to surfaces of type $III_b$. 
\end{teo}

\Proof
$J$ defines a subvariety $X$ in $\PP(1,1,2,3,4,5,6)$.

 We first show
that, for general choice of the entries of $M'$, $X$ has only
rational double points as singularities.

We denote by $V$ the $3-$fold (containing $X$) defined by the $2 \times
2$ minors of the submatrix of $M'$
$$
B=\begin{pmatrix}
z& v & y  & x_1  \\
w& u & P_3 & y  \\
\end{pmatrix}.
$$

Since the condition of having rational double points is well known to be open,
we may assume w.l.o.g.  $P_3=z$. Note that, by row and column
operations on $B$ and analogously on $M'$, one can replace $P_3$ by
a linear combination of $z$ and $x_0^3$ (we are thus working on a subfamily of
codimension at most $1$). 

$V$ is a cone over a quasi smooth scroll in
$\PP(1,2,3,4,5,6)$: it is therefore quasi smooth outside the vertex
$p_0:=(1,0,0,0,0,0,0)$.

$X$ is defined in $V$ by the remaining $3$ equations

\begin{tabular}{ccccc}
$v^2$&=&$z^2P_4$&+&$x_1P_9$\\
$vu$&=&$zwP_4$&+&$yP_9$\\
$u^2$&=&$w^2P_4$&+&$zP_9$.
\end{tabular}

The three above equations describe (as we vary $P_4$ and $P_9$) an open set of a
linear system of divisors on $V$ without base points; by Bertini's theorem,
for general choice of the coefficients of $P_4$ and $P_9$ the surface $X$ is quasi
smooth outside the point
$p_0$. To take care of possible singularities outside $p_0$ we have to intersect
$X$ with the singular locus of $\PP(1,1,2,3,4,5,6)$.

A singular point of $\PP(1,1,2,3,4,5,6)$ has $x_0=x_1=0$: under these
assumptions the equations of $V$ easily force $y=z=0$, and
consequently by the first of the above further three equations we get
the vanishing of the coordinate $v$.

We are left with at most two nonzero coordinates, $u$ and $w$, 
but consider the
last equation ($u^2=\ldots$): for general $P_4$  we
have (up to a rescaling) $P_4=w+\ldots$, thus  we get only one point, exactly the
point
$(0,0,0,0,1,0,1)$. 

This point is in fact a singular point of the ambient space, since the
$\CC^*$ action has in the corresponding point in $\CC^7$ a non
trivial stabilizer $ \cong \ZZ/2\ZZ$. 
Therefore in this point $X$ has an isolated singularity locally
isomorphic to the quotient of a smooth surface by a $\ZZ/2\ZZ$
action: in other words, a singular point of type $A_1$. 

We are left with the vertex $p_0$ of the cone $V$. This point is smooth
for the ambient space $\PP(1,1,2,3,4,5,6)$: we set $\{x_0=1\}$ and work
in  affine coordinates. 

For general choice of the coefficients $P_9$ is invertible at $p_0$: 
since all other entries of $M$ do
vanish in $p_0$ we get that all equations vanish in
$p_0$ and only the pfaffians including $P_9$ are smooth in it.

We get that the Zariski tangent space of $X$ at $p_0$ has dimension
$3$, and more precisely it is $\{x_1=y=z=0\}$.  We eliminate then
(we have set $\{x_0=1\}$) $x_1, y$ and $z$. We obtain the equation
$wv=P_9^{-1}u(u^2-w^2P_4)$: a rational double
point of type $A_2$.

We have  thus shown that, for general choice of the coefficients, $X$ has
only rational double points as singularities.

The projection from $V$ to the quadric cone $\PP(1,1,2)$ given by the
first three variables $x_0,x_1,y$ has  $\PP^1$ as general fibre: for
$x_1 \neq 0$ the equations of $V$ can be explicitly solved, and yield
 $z=y^2/x_1$, $w=y^3/x_1^2$.

The remaining $3$ equations cut clearly two points on the general $\PP^1$
fibre of the projection to $\PP(1,1,2)$: therefore the induced map $X
  \rightarrow \PP(1,1,2)$ has degree $2$. 
 
We have the following recipe to obtain the branch curve:
\begin{itemize}
\item start with the polynomial $z^2P_4+x_1P_9$;
\item substitute $z \mapsto y^2/x_1$; 
\item substitute $w \mapsto y^3/x_1^2$;
\item multiply the result by $x_1^4$.
\end{itemize}
What we get is a polynomial (in fact in the last step we get rid of the
denominators) of degree
$14$ in the variables $x_0,x_1,y$ (remember that we have assumed $v,u$ not to
appear in
$P_9$), therefore  the branch curve is a general curve of degree $14$ in
$\PP(1,1,2)$ contained in the ideal $(y^7, x_1y^6, x_1^2y^4, x_1^3y^3,
x_1^4y^2, x_1^5)$; these curves form a linear system on $\PP(1,1,2)$ with
base point $(1,0,0)$. Its general element has a  $5-$ple point in $(1,0,0)$ and is
smooth elsewhere. 

We blow up the point $(1,0,0)$, take the complete transform ad remove $4$ times
the exceptional divisor. Then we get a triple 
point with tangent cone three times the direction of the exceptional
divisor.

Blowing up again and  removing the new exceptional
divisor twice from the complete transform of the branch curve, we then obtain
a $4$-ple point with (in general) three different tangent directions: after a last
blow up we remain with at most non essential singularities.

Summarizing, our branch curve of degree $14$ has a $5-$tuple point in $(1,0,0)$ with
an infinitely near $(3,3)$-point; by \cite{horIII}, thm. 5.2., the
desingularisation of $X$ is a surface of type
($III_b$). 

In order to count the number of moduli of this subfamily of surfaces of type
(IIIb), we
observe that the two singular points are infinitely near (one condition), the
'infinitely near' triple point is in the direction of the exceptional divisor (one
condition). 

Moreover,
we have the further condition that the coefficient of the monomial
$yx_1^4$ (in the equation of the branch curve) vanishes:  in fact a straightforward
computation shows that the two first conditions force the equation of the branch
curve to be only in the ideal \\
$(y^7, x_1y^6, x_1^2y^4, x_1^3y^3, x_1^4y^2, yx_1^4, x_1^5)$.

%We did note that the curve acquires a double point after the third blow-up:
%this depends on this vanishing (adding this monomial we do not get any
%further singular point).

Under the assumption that $P_3=z$, we have a $38-3=35$
dimensional family of surfaces of type $III_b$; in general we have to write
$P_3=z+\lambda x_0^3$; the projection on $\PP(1,1,2)$
 has clearly still degree $2$ and we have a similar recipe for
 computing the branch curve. We still have a curve as
 described in theorem \ref{IIIb}, still the two base points are
 infinitely near, but the singular point can have multiplicity $4$
 (instead of $5$), so are not in the subfamily and the dimension
 of the whole family is one more, i.e., $36$.

\QED

\begin{oss}
It is clear from the previous proof that the above family 
given by the pfaffians of $M'$ (having dimension $36$)
(cf. \ref{famigliadidimensione36}) is a proper subfamily of the $37$ - dimensional
family of surfaces of type
$III_b$, such that the two essential singularities of the branch curve of
the canonical double cover are infinitely near. (cf. \cite{horIII}). 
We do not have any geometric characterization for this subfamily. 
\end{oss}

\begin{teo}\label{family}
Let $(x_0,x_1,y,z,w,v,u)$ be variables of respective degrees $(1,1,2,3,4,5,6)$, 
Let $M$ be the $6 \times 6$ skew-symmetric matrix
$$M=
\begin{pmatrix}
  0 &t&z& v & y  & x_1  \\
    &0&w& u & P_3& y  \\
    & &0&P_9& u  & v  \\
    & & & 0 &wP_4&zP_4\\
    & & &   & 0  &tP_4\\
-sym& & &   &    & 0
\end{pmatrix}.$$ where the $P_i$'s are homogeneous polynomials of degree $i$
    in the above  variables and $t$ is the parameter on a
    small disk $\Delta_{\varepsilon} \subset \CC$.

For general choice of $P_3,P_4$ and $P_9$ the $4 \times 4$
pfaffians of $M$ define a variety $X \subset \Delta_{\varepsilon} \times
\PP(1,1,2,3,4,5,6)$ whose projection on $\Delta_{\varepsilon}$ is flat, 
with central fibre a surface of type $III_b$ and with general fibre 
a surface of type $II$. 
\end{teo}

\Proof
The flatness of the above family (for general entries) 
follows directly from the flexibility of the format. For $t=0$ the above matrix
equals the matrix $M'$ in thm. \ref{famigliadidimensione36}.

Assume now that $t \neq 0$.

Note that the pfaffians $Pf_{1235}$ and $Pf_{1236}$ are of the form
$tu-\cdots$ and $tv-\cdots$, and that for a general choice of $P_4$, the pfaffian
$Pf_{1256}$ can be written as $t^2w-...$. 
Therefore, for $t \neq 0$, we can eliminate the variables $u,v,w$,
and $R \cong \CC[x_0,x_1,y,z]/J$ for a suitable ideal $J$; a straightforward 
calculation shows that $J$ is a principal ideal generated by the
equation obtained from $Pf_{1345}$ after eliminating the variables $u,v,w$ using 
$Pf_{1235}$, $Pf_{1236}$ and $Pf_{1256}$.  

This is a polynomial of degree $9$, so the surface is birational to a
hypersurface of degree $9$ in $\PP(1,1,2,3)$, 
 whence we obtain a surface of type $II$. 

\QED

Our main theorem (\ref{Degree6}) follows right away from the above.

{\em Authors' addresses} \\

\noindent
PD. Dr. I. Bauer, Prof. Dr. F. Catanese \\
Lehrstuhl Mathematik VIII \\
Universit\"at Bayreuth \\
Universit\"atsstr. 30, D-95447 Bayreuth\\

\noindent
Dr. R. Pignatelli \\
Dipartimento di Matematica \\
Via Sommarive, 14, I- 38050 Povo (TN)
Universita'di Trento
\end{document}